\documentclass[12pt]{amsart}
\usepackage{amssymb}

\newtheorem{Thm}{Theorem}[section]

\newtheorem{Cor}[Thm]{Corollary}
\newtheorem{Lem}[Thm]{Lemma}
\newtheorem{Prop}[Thm]{Proposition}

\theoremstyle{definition}

\theoremstyle{remark}

\def \cal{\mathcal}

\def\Cdb{\mathbb C}

\def\Rdb{\mathbb R}
\def\Tdb{\mathbb T}

\begin{document}

\title{A solution to the problem of $L^p-$maximal regularity}
\author{N. J. Kalton}
\address{Department of Mathematics \\
University of Missouri-Columbia \\
Columbia, MO 65211
}
\thanks{The first author was partially supported by NSF grant
DMS-9870027}

\email[N. J. Kalton]{nigel@math.missouri.edu}

\author{G. Lancien}
\address{Equipe de Math\'ematiques - UMR 6623, Universit\'e de
Franche-Comt\'e, F-25030 Besan\c con cedex}

\email{GLancien@vega.univ-fcomte.fr}

\subjclass{Primary:47D06, Secondary: 46B03, 46B52, 46C15}

\begin{abstract}
We give a negative solution to the problem of the $L^p$-maximal regularity on
various classes of Banach spaces including $L^q$-spaces with $1<q \neq
2<+\infty$.
\end{abstract}

\maketitle

\section{Introduction}

In this paper we consider the following abstract Cauchy
problem:
$$\left\{
\begin{array}{ll}
u'(t)+B(u(t))=f(t)\ \ \ \ {\rm for\ } 0\leq t<T \\
u(0)=0
\end{array}\right.$$
where $T \in (0,+\infty)$, $-B$ is the infinitesimal generator of a bounded
analytic semigroup on a complex Banach space $X$ and $u$ and $f$ are $X$-valued
functions on $[0,T)$. Suppose $1<p<\infty.$ $B$ is said to satisfy {\it
$L^p-$maximal
regularity} if whenever $f\in L^p([0,T);X)$ then the solution
$$ u(t) =\int_0^t e^{-(t-s)B}f(s)\,ds$$
satisfies  $u'\in L^p([0,T);X).$  It is known that $B$ has $L^p-$ maximal
regularity for some $1<p<\infty$ if and only if it has $L^p-$maximal
regularity for every $1<p<\infty$ \cite{DE}, \cite{DO}, \cite{SO}.  We
thus say simply that $B$ satisfies {\it maximal regularity (MR).}

The question whether $B$ satisfies maximal regularity
has been extensively studied.  De Simon
\cite{DE} showed that $B$ always satisfies (MR) if $X$
is a Hilbert space. We also mention the early work of
 Grisvard \cite{G} using
interpolation spaces between $X$ and the domain $D(B).$ More recently,
  Dore and  Venni \cite{DV} showed that $B$ satisfies (MR) if
$X$ is an
UMD Banach space and $B$ admits bounded imaginary powers with an estimate
$||B^{is}||\leq Ke^{\theta |s|}$ for some $0\leq
\theta<\pi/2.$
 The question we address was first asked by H. Br\'ezis (see
\cite{COU}) and is the following: under which conditions on the Banach space
$X$ does every negative generator of a bounded analytic semigroup on $X$
satisfies (MR).
  Let us say that $X$ has the {\it the maximal regularity property (MRP)}
if $B$ satisfies (MR)  whenever $-B$ is the generator of a bounded
analytic
semigroup.  The result of De Simon cited above implies that Hilbert
spaces satisfy (MRP).  Note also that by a theorem of Lotz \cite{LO}
every strongly continuous semigroup on $L^{\infty}$ is uniformly
continuous and this implies $L^{\infty}$ has (MRP).
 On the other hand,
Coulhon and Lamberton \cite{COU} exhibited counterexamples on
$X=L^p(\Rdb;E)$ whenever $1<p<\infty$ and
$E$ is not an UMD Banach space. More recently,  Le Merdy \cite{LE1} found
counterexamples on other fundamental spaces such as $L^1(\Tdb)$,
$C(\Tdb)$ and
${\mathcal K}(\ell_2)$. For several years it has been an open question
whether the spaces $L^p$ have (MRP) or even whether every UMD-space has
(MRP).

In this paper,
we provide a fairly complete answer to
this question. In fact we will show that (MRP),
 up to isomorphism, characterizes Hilbert spaces
among spaces with an unconditional basis or (more generally)
 separable Banach lattices. We also extend the Coulhon-Lamberton result
cited above by  showing that $L^p(\Rdb;E)$ for $1\le p<\infty$ can never
have
(MRP) unless $p=2$ and $E$ is a Hilbert space.

We would like to
mention that our constructions have been initially inspired by a very useful
transference principle for maximal regularity proved by C. Le Merdy \cite{LE1},
although we only use a very simple version of it (see our Proposition
\ref{trans}).

This work was done during a visit of the  first author to the Department
of Mathematics at
 the Universit\'e de
Franche-Comte in spring 1999; he would like to thank the Department for
its warm hospitality.

\section{Notation and basic facts}

We will mainly adopt the notation introduced in \cite{LE1}.

Let $X$ be a complex Banach space, $1<p<\infty$ and $0<T<\infty$. We denote by
${\cal A}_T$ the differentiation $d/dt$ on
$L^p([0,T);X)$ with domain
$D({\cal A}_T)=W_0^{1,p}([0,T);X)=\{u \in W^{1,p}([0,T);X):\ u(0)=0\}$. Let now
$-B$ be the generator of a bounded analytic semigroup on $X$. The operator
$I_{L^p}\otimes B$ defined on
$L^p([0,T)) \otimes D(B)$ is closable and we denote by ${\cal B}$ its closure.
${\cal B}$ can also be described by $D({\cal B})=L^p([0,T);D(B))$ and $({\cal
B}u)(t)= B(u(t))$ for $u$ in $D({\cal B})$. We say that $B$ satisfies
$L^p$-maximal regularity (on $[0,T)$) if the operator ${\cal A}_T + {\cal
B}$
with domain $D({\cal A}_T) \cap D({\cal B})$ has a bounded inverse; this
is equivalent to the formulation in the introduction, and  obviously
independent of $T.$ By a result in
\cite{DA} $L^p-$maximal regularity is also equivalent to the statement
that ${\mathcal A}_T+{\mathcal B}$ is closed, and also to an inequality
of the form
$$ \|{\mathcal A}_Tu\| \le C\|{\mathcal A}_Tu+{\mathcal B}u\|$$ for $u\in
D({\mathcal A}_T)\cap D(\mathcal B).$
 As remarked in the introduction, it is
known (see \cite{DE}, \cite{DO}, \cite{SO}) that this property does not depend
on
$1<p<\infty.$
 Thus we will simply say that $B$ satisfies (MR) and work
only on
$L^2([0,T);X)$. Then we say that a Banach space $X$ has the {\it maximal
regularity property (MRP)}, if every $B$ such that $-B$
generates a bounded analytic semigroup on $X$ satisfies (MR).

 We will also use the following terminology. A closed densely defined operator
$B$ on a Banach space $X$ is said to be sectorial of type $\omega$, where
$0<\omega<\pi$, if the spectrum $\sigma(B)$ of $B$ is
included in $\overline\Sigma_\omega$, where $\Sigma_\omega=\{z\in \Cdb:\ |{\rm
Arg}(z)|<\omega \}$ and for every $\omega<\theta<\pi$ there exists
$C_{\theta}>0$ so that for any $\lambda \notin \Sigma_{\theta}$ we have
$\|(\lambda-B)^{-1}\|\le C_{\theta}|\lambda|^{-1}.$
Notice that $-B$ generates a bounded analytic semigroup on $X$ if and only if
$B$ is sectorial of type $\omega$, for some $\omega <
\pi/2$ (see \cite{TA} or \cite{LE2} for details).

Our next result can be regarded as a form of transference to the
circle. For
$f\in
L^1([0,2\pi);X)$ we define the Fourier coefficients $\hat f(n)$ in the
usual way
$$ \hat f(n)= \frac{1}{2\pi}\int_0^{2\pi} f(t)e^{-int}dt$$
for $n\in \mathbb Z.$

\begin{Prop}\label{trans} Let $X$ be a Banach space and let $-B$ be an
invertible generator of a bounded analytic semigroup on $X$. Assume that $B$
satisfies (MR).   Then there is a constant $C$ so that for any $X$-valued
trigonometric polynomial $f(t)=\sum_{k=-N}^N \hat f(n) e^{int}$ we have
$$ \left(\int_0^{2\pi}\| \sum_{n\in \mathbb Z} in (in +B)^{-1} \hat f(n)
e^{int}\|^2\frac{dt}{2\pi}\right)^{1/2} \le
C\left(\int_0^{2\pi}\|
f(t)\|^2\frac{dt}{2\pi}\right)^{1/2}.$$
\end{Prop}

\begin{proof} Denote by $(e^{-tB})_{t>0}$ the semigroup generated by $-B$.
For any trigonometric polynomial $f$ define $g\in L^2([-2\pi,2\pi);X)$
by:
$$\left\{
\begin{array}{ll}
g(s)= (1-e^{-2\pi B})^{-1}f(s+2\pi)\ \ \ \ {\rm for\ } -2\pi\leq s\leq 0
\\ g(s)=f(s)\ \ \ \ {\rm for\ } 0< s\leq 2\pi
\end{array}\right.$$
Since $B$ has a negative exponential type, the definition of $g$ makes sense.
Moreover there is a constant $C>0$ independent of $f$ so that
$$||g||_{L^2([-2\pi,2\pi);X)} \leq
C||f||_{L^2([0,2\pi);X)}.$$
Now we solve
$$\left\{
\begin{array}{ll}
u'(t)+B(u(t))=g(t)\ \ \ \ {\rm for\ } -2\pi\leq t<2\pi \\
u(-2\pi)=0
\end{array}\right.$$
Elementary calculations show that $u'=\sum_{n\in \mathbb Z} in
(in+B)^{-1}\hat f(n) e^{int}$ on
$(0,2\pi)$. Now, the fact that
$B$ has (MR) yields the result.
\end{proof}

Next, is an elementary lemma about (MRP) that we will use extensively:

\begin{Lem}\label{complemented} Let $X$ be a Banach space and $Y$ be a
complemented subspace of
$X$. Assume that $X$ has the (MRP) then $Y$ has the (MRP).
\end{Lem}

\begin{proof} Assume $X=Y \oplus Z$ and that $B$ is a sectorial operator of
type $<\pi/2$ on $Y$ which fails (MR).
 Then the operator $B'$, defined on $X$ by $D(B')=
D(B) \oplus Z$ and $B'x=By$ when  $x=y+z$ with $(y,z) \in D(B) \times
Z$
 provides a counterexample to (MR) on $X$.
\end{proof}

The operators that we will use will be multipliers associated with various
Schauder decompositions.  Let us introduce some notation for that
purpose.
If $F \subset X$, we denote by $[F]$ the closed linear span of $F$. Let
$(E_n)_{n\geq1}$ be a sequence of closed subspaces of $X$. Assume that
$(E_n)_{n\geq1}$ is a Schauder decomposition of $X$ and let $(P_n)_{n\geq1}$ be
the associated sequence of projections from $X$ onto $E_n$. For
convenience  we
will also denote this Schauder decomposition by $(E_n,P_n)_{n\geq1}$.
Notice that the spaces $Z_n=P_n^*(X^*)$ form a Schauder decomposition of
the subspace $Z=[\cup_{n=1}^{\infty} Z_n]$ of $X^*.$

 Let now
$(b_n)_{n\geq1}$ be a sequence of complex numbers. We define the
(possibly unbounded) operator
$M(b_n;E_n)$ with domain
$D(M(b_n;E_n))= \{x\in X~{\rm such~that}~
 \sum b_nP_nx~{\rm converges~in~}X\}$ by
 $M(b_n;E_n)x= \sum
b_nP_nx$.

The following lemma is elementary (see \cite{BC} or \cite{V} for a proof in the
case of a Schauder basis).

\begin{Lem}\label{mult}

\noindent (i) $M(b_n;E_n)$ is a closed densely defined operator.

\noindent (ii) If  $b_1 >0$ and $(b_n)$ is an increasing sequence of reals,
 then $M(b_n;E_n)$ is invertible and sectorial of type
$\omega$ for any
$\omega \in (0,\pi)$.
\end{Lem}

\section{The results}

We first establish a necessary condition for spaces with a Schauder
decomposition to have the maximal regularity property.

\begin{Thm}\label{rad} Let $(E_n,P_n)_{n \geq 1}$ be a Schauder
decomposition of the Banach space $X$. Let $Z_n=P_n^*X^*$ and
$Z=[\cup_{n=1}^{\infty}Z_n].$   Assume $X$ has (MRP). Then there is a
constant
$C>0$ so that
whenever $(u_n)_{n=1}^N$ are such that $u_n\in [E_{2n-1},E_{2n}]$ and
$(u_n^*)_{n=1}^N$ are such that $u_n^*\in [Z_{2n-1},Z_{2n}]$ then
$$ \left(\int_0^{2\pi}\|\sum_{n=1}^NP_{2n}u_n
e^{i2^nt}\|^2\frac{dt}{2\pi}\right)^{1/2}\le C\left( \int_0^{2\pi} \|
\sum_{n=1}^Nu_ne^{i2^nt}\|^2\frac{dt}{2\pi}\right)^{1/2}$$
and
$$ \left(\int_0^{2\pi}\|\sum_{n=1}^NP^*_{2n}u^*_n
e^{i2^nt}\|^2\frac{dt}{2\pi}\right)^{1/2}\le C\left( \int_0^{2\pi} \|
\sum_{n=1}^Nu^*_ne^{i2^nt}\|^2\frac{dt}{2\pi}\right)^{1/2}.$$
\end{Thm}

\begin{proof} Let $(a_n)_{n\geq1}$ and  $(b_n)_{n\geq1}$ be two
sequences defined by
$$a_{2n-1}=b_{2n-1}=b_{2n}=2^{n-1}~Ê~Ê {\rm and}~Ê~Ê a_{2n}=2^n.$$
We let $A=M(a_n;E_n)$ and $B=M(b_n;E_n).$  It is easy to see that
$$ (i2^n+A)^{-1}u_n = (i2^n+2^{n-1})^{-1}P_{2n-1}u_n +
(i2^n+2^n)^{-1}P_{2n}u_n$$ while
$$ (i2^n+B)^{-1}u_n= (i2^n+2^{n-1})^{-1}u_n.$$
Hence
$$ i2^n (i2^n+B)^{-1}u_n-i2^n(i2^n+A)^{-1}u_n =
\frac{i}{(i+1)(2i+1)}P_{2n}u_n.$$

If we assume that $X$ has (MRP) then both $A$ and $B$ satisfy (MR) and so
we can apply Proposition \ref{trans} to each in turn for the polynomial
$f(t)=\sum_{n=1}^N u_n e^{i2^nt}.$  Subtracting gives us the first
estimate.

The second estimate follows by duality.  More precisely the operators
$f\mapsto \sum_{n\in\mathbb Z}in (in +A)^{-1} \hat f(n)e^{int}$ and
$f\mapsto
\sum_{n\in\mathbb Z} in (in +B)^{-1} \hat f(n) e^{int}$ can be extended
to bounded linear operators on $L^2([0,2\pi);X).$ Taking
adjoints and restricting to the subspace $L^2([0,2\pi);Z)$ one easily
obtains similar estimates in the dual. \end{proof}

We first examine two important examples.

\begin{Cor}\label{czero}
$c_0$ and $\ell_1$ fail the (MRP).
\end{Cor}

\begin{proof} Denote  by $(x_n)_{n \geq 1}$ the canonical basis of $c_0$ and
let $s_n=x_1+..+x_n$. $(s_n)_{n \geq 1}$ is a Schauder basis of $c_0$ which is
usually called the summing basis of $c_0$. We now apply Theorem \ref{rad}
with the sequence of projections $(P_n)$ associated with the Schauder basis
$(s_n)$ and $u_n=s_{2n}-s_{2n-1}$. Then we obtain that there is $C>0$ so that for every
$N\ge 1,$:
$$\left(\int_0^{2\pi} \|\sum_{n=1}^N
s_{2n}e^{i2^nt}\|^2\frac{dt}{2\pi}\right)^{1/2} \le
C\left(\int_0^{2\pi}\|\sum_{n=1}^N (s_{2n}-s_{2n-1})e^{i2^nt}
||^2\frac{dt}{2\pi}\right)^{1/2}.$$
The right-hand side is equal to $C$ but, considering only the
first
co-ordinate of the left-hand side with respect to the canonical basis, we
have
$$\left(\int_0^{2\pi} \|\sum_{n=1}^N
s_{2n}e^{i2^nt}\|^2\frac{dt}{2\pi}\right)^{1/2} \ge N^{1/2}.$$
This is a contradiction.

\smallskip
Assume now that $\ell_1$ has the (MRP). Let $(v_n)_{n \geq 1}$ be the
coordinate functionals associated with the summing basis $(s_n)$ of $c_0$. The
closed linear space $Y$ spanned in $\ell_1$ by the sequence $(v_n)$ is of
codimension 1 in $\ell_1$ and hence is isomorphic to $\ell_1.$   The
bi-orthogonal functionals $(v_n^*)$ in $Y^*$ are equivalent to the
summing basis of $c_0.$  Hence using the same calculation as above and
the second inequality of Theorem \ref{rad} we again get a contradiction.
\end{proof}

We now explain the consequences of Theorem \ref{rad} when $X$ admits
an unconditional basis.

\begin{Thm}\label{uncbasis} A Banach space with an unconditional basis has the
(MRP) if and only if it is isomorphic to $\ell_2$.
\end{Thm}

\begin{proof} The idea is to show that if $X$ has the (MRP) and an
unconditional basis $(x_n)_{n\geq 1}$, then for every  permutation $\pi$ of the
integers and for every block basis
$(u_j)_{j
\geq 1}$  of $(x_{\pi(n)})_{n \geq 1}$ the closed subspace of $X$
spanned by
the $u_j$'s is complemented in $X$. Once we have shown this the proof is
completed by using a theorem of Lindenstrauss
and  Tzafriri (\cite{LT}, see also \cite{LIT1} Theorem 2.a.10) which
asserts that
$(x_n)_{n
\geq 1}$ must be equivalent to the canonical basis of $c_0$ or $\ell_p$ for
some $p$ in $[1,\infty)$. Then, by Corollary \ref{czero}, $(e_n)$ is
equivalent to  the canonical basis of $\ell_p$ for
some $p$ in $(1,\infty)$. Now, if
$1<p\neq 2<\infty$,
$\ell_p$ admits an unconditional basis which is not equivalent to any of the
canonical bases of
$c_0$ or $\ell_q$ where $1\leq q<\infty$. Indeed  Pe\l czy\'nski
\cite{PE} showed
that, for $1<p<\infty$, $\ell_p$ is isomorphic to $(\displaystyle\sum_{n \geq
1}\oplus \ell_2^n)_p$.

So assume, as we may, that $(x_n)_{n \geq 1}$ is a normalized
1-unconditional
basis of $X$ and that $(u_n)_{n \geq 1}$ is a normalized block basis of
$(x_n)_{n \geq 1}$, with
$$\forall n \geq 1,\  u_n= \sum_{r_n+1}^{r_{n+1}} a_je_j,$$
where $1=r_1<r_2<..<r_n<r_{n+1}<..$ and $(a_j)_{j\geq 1} \subset \Cdb$. For $n
\geq 1$, let
$X_n= [x_{r_n+1},..,x_{r_{n+1}}]$ and $E_{2n}=[u_n]$. Then $(X_n)$ is an
unconditional Schauder decomposition of $X$ with associated projections
$P_n,$ say. Now, by the Hahn-Banach theorem there is a norm-one
projection
$R_n:X_n\to E_{2n}.$ Let $E_{2n-1}=R_n^{-1}(0).$  Then $(E_n)$ is a
Schauder decomposition of $X$ with associated projections
$Q_{2n-1}=(I-R_n)P_n$ and $Q_{2n}=R_nP_n.$  We now apply Theorem
\ref{rad} and exploit the unconditionality of the Schauder decomposition
$(X_n).$
 There is a constant
$C$ so that if
$y$ is in the linear span of the $(x_n)_{n\ge 1}$ then
$$ \|\sum_{n=1}^{\infty}Q_{2n}y\| \le C\|y\|.$$  This implies that
$[u_n]_{n\ge 1}$ is complemented in $X.$

Clearly the same reasoning can be applied to any permutation of the basis
$(x_n)$ so that the proof is complete.
\end{proof}

Although this will be included in  further and more general statements let us
point out that this already solves our problem for the spaces
$L^p(0,1)$:

\begin{Cor}\label{Lp} Let $1 \leq p\leq\infty$. Then $L^p(0,1)$ has the (MRP)
if and only if $p=2$ or $p=\infty$.
\end{Cor}

\begin{proof} For $1<p<\infty$, the Haar system is known to be an
unconditional basis of $L^p(0,1)$ (\cite{PA}, see also \cite{LIT2}). So the
result follows from the preceding Theorem. The fact that $L^1$ fails (MRP)
was proved by C. Le Merdy in \cite{LE1}. Notice that $L^1$ contains a
complemented copy of $\ell_1$, so this result can be derived from Lemma
\ref{complemented} and Corollary \ref{czero}.
\end{proof}

We now extend Theorem \ref{uncbasis} to the case of a space with an
unconditional Schauder decomposition.

\begin{Thm}\label{unc} Let $X$ be a Banach space with an unconditional
decomposition
$(F_n,P_n)_{n \geq 1}$. Assume that $X$ has the (MRP).

\noindent Then $X$ is isomorphic
to
$(\displaystyle\sum_{n=1}^\infty \oplus F_n)_{\ell_2}$.
\end{Thm}

\begin{proof}
It suffices to show that if $u_n\in F_n$ with $\|u_n\|=1$ then $\sum
a_nu_n$ converges if and only if $\sum |a_n|^2<\infty.$   As above in
the proof of Theorem \ref{unc}, let
$R_n$ be a norm-one projection of $F_n$ onto $[u_n].$   Then let
$E_{2n}=[u_n]$ and $E_{2n-1}=R_n^{-1}(0).$  Reasoning exactly as in
Theorem \ref{unc} gives that $[u_n]_{n\ge 1}$ is complemented in $X.$
But this subspace has an unconditional basis and so Theorem \ref{unc}
yields that $(u_n)$ is equivalent to the canonical basis of $\ell_2.$
\end{proof}

Our next theorem completes the important counterexamples obtained by
Coulhon and  Lamberton \cite{COU}.

\begin{Thm}\label{bochner} Suppose $X$ is a Banach space and $1\le
p<\infty.$ Then the Banach space
$L^p((0,1);X)$ has the
(MRP)  if and only if $p=2$ and $X$ is isomorphic to a Hilbert space.
\end{Thm}

\begin{proof} We first note that $L^p[0,1]$ is complemented in
$L^p((0,1);X)$ so that if the latter has (MRP) then $p=2$ by
Corollary \ref{Lp}.  Assume that
$L^2((0,1);X)$ has the
(MRP); we will show that $X$ is isomorphic to a Hilbert space
(the opposite implication is
due to  de Simon \cite{DE}). By \cite{COU}, $X$ must have the UMD property.
In particular, $X$ does not contain the $\ell_1^n$'s uniformly. Hence by
Pisier's theorem \cite{PI2}, the space ${\rm
Rad}(X)=\overline{[\varepsilon_n]_{n\geq1} \otimes X}$ is complemented in
$L^2((0,1);X)$ (here $\varepsilon_n$ is a standard Rademacher function).
Therefore, by Lemma
\ref{complemented},
${\rm Rad}(X)$ has the
(MRP). Now,
$(E_n)_{n\geq1}=(\varepsilon_n \otimes X)_{n\geq1}$ is an unconditional
Schauder  decomposition of ${\rm Rad}(X)$. So, by Theorem \ref{unc}, ${\rm
Rad}(X)$ must be isomorphic to $(\sum\oplus \,(\varepsilon_n \otimes
X))_{\ell_2}$. Finally, it follows from Kwapien's theorem \cite{K} that $X$ is
isomorphic to a Hilbert space.
\end{proof}

We now extend Theorem \ref{uncbasis} and Corollary \ref{Lp} to the setting of
Banach lattices. All the notions on Banach lattices that we will use can be
found in
\cite{LIT2} Chapters 1.a and 1.b.

\begin{Thm}\label{lat} An order continuous Banach lattice has the (MRP) if and
only if it is isomorphic to a Hilbert space.
\end{Thm}

\begin{proof} Let $X$ be an order continuous Banach lattice. By \cite{LIT2},
Lemma 1.b.13, it is enough to show that every normalized sequence of disjoint
elements of $X$ is equivalent to the canonical basis of $\ell_2$. So let
$(f_n)_{n\geq 1}$ be such a sequence in $X$. Then $X$ admits an unconditional
Schauder decomposition $(E_n)_{n\geq1}$ such that the $E_n$'s are ideals of
$X$ and for all $n\geq 1$, $f_n\in E_n$. Now, by  Theorem \ref{unc} $X$ is
isomorphic to $(\sum\oplus E_n )_{\ell_2}$ and $(f_n)$ is equivalent to the
canonical basis of $\ell_2.$
\end{proof}

\begin{Cor}\label{seplat} A separable Banach lattice has the (MRP) if and only
if it is isomorphic to a Hilbert space.
\end{Cor}

\begin{proof} Let $X$ be a separable Banach lattice which is not order
continuous. Then $X$ is not $\sigma$-complete (see \cite{LIT2} Proposition
1.a.7) and by a result of P. Meyer-Nieberg (\cite{MN}, see also \cite{LIT2}
Theorem 1.a.5) $X$ contains a subspace isomorphic to $c_0$. Since $X$ is
separable, it follows from Sobczyk's Theorem \cite{S} that this subspace is
complemented in
$X$. So, by Lemma \ref{complemented} and Corollary \ref{czero}, $X$ does not
have the (MRP). Then, the preceding Theorem concludes our proof.
\end{proof}

\section {Final remarks}

\medskip\noindent
1) One can also consider the problem of the $L^p$-maximal regularity on the
half line $[0,+\infty)$, which has in general a different answer (see
\cite{LE1} for an example). But it follows from Theorem 2.4. in \cite{DO} that
all our results remain valid in this slightly different setting.

\smallskip\noindent
2) We do not know if there is a non-Hilbertian subspace of an
$L^p$-space
($1\leq p<\infty$) with  (MRP).

\smallskip\noindent
3) We do not know if every space with a basis and (MRP) is isomorphic to
a Hilbert space.


\vskip 2cm

\begin{thebibliography}{LIT1}

\bibitem{BC} J.B. Baillon and P. Cl\'ement, Examples of unbounded
imaginary powers of operators, J. Funct. Anal. 100 (1991), 419-434.


\bibitem{COU} T. Coulhon and D. Lamberton, R\'egularit\'e $L^p$ pour les
\'equations d'\'evolution, S\'eminaire d'Analyse Fonctionnelle Paris VI-VII
(1984-85),155-165.

\bibitem{DA} G. Da Prato and P. Grisvard, Sommes d'op\'erateurs
lin\'eaires et
\'equations diff\'erentielles op\'erationnelles, J. Math. Pures Appl. 54
(1975), 305-387.

\bibitem{DE} L. De Simon, Un' applicazione della theoria degli integrali
singolari allo studio delle equazioni differenziali lineare astratte del primo
ordine, Rend. Sem. Mat., Univ. Padova (1964), 205-223.

\bibitem{DO} G. Dore, $L^p$ regularity for abstract differential equations
(In ``Functional Analysis and related topics", editor: H. Komatsu), Lect.
Notes in Math. 1540, Springer Verlag (1993).

\bibitem{DV} G. Dore, A. Venni, On the closedness of the sum of two closed
operators, Math. Z. 196 (1987), 189-201.

\bibitem{G} P. Grisvard, Equations diff\'erentielles abstraites, Ann. Sci.
Ecole Norm. Sup. (4) 2 (1969), 311-395.

\bibitem{K} S. Kwapien, Isomorphic characterizations of inner product spaces
by orthogonal series with vector coefficients, Studia Math. 44 (1972), 583-595.
\bibitem{LE1} C. Le Merdy, Counterexamples on $L^p$-maximal regularity, Math.
Z. 230 (1999), 47-62.

\bibitem{LE2} C. Le Merdy, $H^\infty$-functional calculus and applications to
maximal regularity, Publications Math\'ematiques de Besan\c con, Fasicule 16
(1999).

\bibitem{LT} J. Lindenstrauss and L. Tzafriri, On the complemented
subspaces problem, Israel J. Math. 9 (1971), 263-269.

\bibitem{LIT1} J. Lindenstauss and L. Tzafriri, {\it Classical Banach
spaces I}, Springer-Berlin (1977).

\bibitem{LIT2} J. Lindenstauss and L. Tzafriri, {\it Classical Banach
spaces II}, Springer-Berlin (1979).

\bibitem{LO} H.P. Lotz, Uniform convergence of operators on $L^\infty$ and
similar spaces, Math. Z. 190 (1985), 207-220.

\bibitem{MN} P. Meyer-Nieberg, Charakterisierung einiger topologischer und
ordnungstheoritischer Eigenschaften von Banachverbanden mit Hilfe disjunkte
Folgen, Arch. Math. 24 (1973), 640-647.

\bibitem{PA} R.E. Paley, A remarkable sequence of orthogonal functions,
Proc. London Math. Soc. 34 (1932) 241-264.


\bibitem{PE} A. Pe\l czy\'nski, Projections in certain Banach spaces, Studia Math.
19 (1960), 209-228.


\bibitem{PI2} G. Pisier, Holomorphic semigroups and the geometry of Banach
spaces, Annals of Math. (2) 115 (1982), 375-392.

\bibitem{S} A. Sobczyk, Projection of the space $m$ on its subspace $c_0$,
Bull. Amer. Math. Soc. 47 (1941), 938-947.

\bibitem{SO} P.E. Sobolevskii, Coerciveness inequalities for abstract parabolic
equations (translations), Soviet Math. Dokl. 5 (1964), 894-897.

\bibitem{TA} H. Tanabe,{ \it Equations of evolution,} Pitman, London-San
Francisco-Melbourne, 1979.

\bibitem{V} A. Venni, A counterexample concerning imaginary powers of linear
operators, (In ``Functional Analysis and related topics", editor: H. Komatsu),
Lect. Notes in Math. 1540, Springer Verlag (1993).




\end{thebibliography}
\end{document}